\documentclass[11 pt]{article}
\usepackage[T1]{fontenc}
\usepackage{stackengine, IEEEtrantools}
\usepackage{cite, authblk}
\usepackage{amsmath,amssymb,amsfonts}
\usepackage{algorithmic}
\usepackage{graphicx}
\usepackage{textcomp}
\usepackage{xcolor}
\usepackage{algorithm,algorithmic}
\usepackage{mathtools, stackengine, mathrsfs, setspace, bm, multirow, tikz, pgfplots, adjustbox, calc, url}
\usepackage{xparse}
\pgfplotsset{width=8.1cm, legend style={font=\small}}
\newsavebox{\fminipagebox}
\NewDocumentEnvironment{fminipage}{m O{\fboxsep}}
 {\par\kern#2\noindent\begin{lrbox}{\fminipagebox}
  \begin{minipage}{#1}\ignorespaces}
 {\end{minipage}\end{lrbox}%
  \makebox[#1]{%
    \kern\dimexpr-\fboxsep-\fboxrule\relax
    \fbox{\usebox{\fminipagebox}}%
    \kern\dimexpr-\fboxsep-\fboxrule\relax
  }\par\kern#2
 }
\begingroup\edef\x{\endgroup
  \mathchardef\mathdollar=\the\numexpr"7000+\the\mathdollar\relax
}\x
\DeclareMathAlphabet{\mathit}{T1}{cmr}{m}{it}

\makeatletter
\let\MYcaption\@makecaption
\makeatother
\usepackage{caption}
\DeclareCaptionLabelSeparator{colonquad}{:\quad}
\DeclareCaptionStyle{fig-top}% 
	[justification=raggedright,indention=2cm, font=footnotesize, labelsep=colonquad]{}
\DeclareCaptionStyle{fig}% 
	[justification=raggedright,indention=2cm, font=footnotesize, labelsep=colonquad]{}	
\makeatletter
\let\@makecaption\MYcaption
\makeatother

\pgfplotsset{
    every first x axis line/.style={},
    every first y axis line/.style={},
    every first z axis line/.style={},
    every second x axis line/.style={},
    every second y axis line/.style={},
    every second z axis line/.style={},
    first x axis line style/.style={/pgfplots/every first x axis line/.append style={#1}},
    first y axis line style/.style={/pgfplots/every first y axis line/.append style={#1}},
    first z axis line style/.style={/pgfplots/every first z axis line/.append style={#1}},
    second x axis line style/.style={/pgfplots/every second x axis line/.append style={#1}},
    second y axis line style/.style={/pgfplots/every second y axis line/.append style={#1}},
    second z axis line style/.style={/pgfplots/every second z axis line/.append style={#1}}
}

\makeatletter
\def\pgfplots@drawaxis@outerlines@separate@onorientedsurf#1#2{%
    \if2\csname pgfplots@#1axislinesnum\endcsname
    \else
    \scope[/pgfplots/every outer #1 axis line,
        #1discont,decoration={pre length=\csname #1disstart\endcsname, post length=\csname #1disend\endcsname}]
        \pgfplots@ifaxisline@B@onorientedsurf@should@be@drawn{0}{%
            \draw [/pgfplots/every first #1 axis line] decorate {
                \pgfextra
                \pgfplotspointonorientedsurfaceabsetupfor{#2}{#1}{\pgfplotspointonorientedsurfaceN}%
                \pgfplots@drawgridlines@onorientedsurf@fromto{\csname pgfplots@#2min\endcsname}%
                \endpgfextra 
                };
        }{}%
        \pgfplots@ifaxisline@B@onorientedsurf@should@be@drawn{1}{%
            \draw [/pgfplots/every second #1 axis line] decorate {
                \pgfextra
                \pgfplotspointonorientedsurfaceabsetupfor{#2}{#1}{\pgfplotspointonorientedsurfaceN}%
                \pgfplots@drawgridlines@onorientedsurf@fromto{\csname pgfplots@#2max\endcsname}%
                \endpgfextra 
                };
        }{}%
    \endscope
    \fi
}%
\makeatother
\newcommand\blfootnote[1]{%
  \begingroup
  \renewcommand\thefootnote{}\footnote{#1}%
  \addtocounter{footnote}{-1}%
  \endgroup
}
\def\BibTeX{{\rm B\kern-.05em{\sc i\kern-.025em b}\kern-.08em
    T\kern-.1667em\lower.7ex\hbox{E}\kern-.125emX}}

\allowdisplaybreaks

\title{A Distributionally Robust Optimization Approach\\ for Unit Commitment in Microgrids
}
\author{Ogun Yurdakul}
\author{Fikret Sivrikaya}
\author{Sahin Albayrak}
\affil{Department of Electrical Engineering and Computer Science\\
Technical University of Berlin, Berlin, Germany\\
Email: \{yurdakul, fikret.sivrikaya, sahin.albayrak\}@tu-berlin.de}
\date{}                     %% if you don't need date to appear
\begin{document}
\maketitle
\blfootnote{This work was supported in part by the Research Council of Norway under the ``LUCS'' project, and by the German Federal Ministry for Economic Affairs and Energy under Grant 03EI6004B.}
\begin{abstract}
This paper proposes a distributionally robust unit commitment approach for microgrids under net load and electricity market price uncertainty. The key thrust of the proposed approach is to leverage the Kullback-Leibler divergence to construct an ambiguity set of probability distributions and formulate an optimization problem that minimizes the expected cost brought about by the worst-case distribution in the ambiguity set. The proposed approach effectively exploits historical data and capitalizes on the k-means clustering algorithm---in conjunction with the soft dynamic time warping score---to form the nominal probability distribution and its associated support. A two-level decomposition method is developed to enable the efficient solution of the devised problem. We carry out representative studies and quantify the relative merits of the proposed approach vis-\`a-vis a stochastic optimization-based model under different divergence tolerance values. \\\\
\end{abstract}
\section{Introduction}\label{1}
A microgrid is a cluster of loads, thermal generation resources (\textit{TGR}s), variable energy resources (\textit{VER}s), and electric storage resources that operate in coordination to supply electricity in a reliable manner. Typically integrated to its host power system at the distribution level, a microgrid is---for all intents and purposes---a microcosm of a bulk power system that retains most of its innate operational characteristics.\par
Similar to bulk power systems, the short-term planning of microgrids can be determined via unit commitment (\textit{UC}) and economic dispatch (\textit{ED}) decisions \cite{drogmb:79}. The \textit{UC} problem seeks minimum cost strategies to determine the commitment statuses of \textit{TGR}s based on expected load, equipment limitations, and operational policies. The equipment limitations of \textit{TGR}s and the inter-temporal constraints of microgrid physical asset operations render \textit{UC} a time-coupled problem and necessitate that the \textit{UC} decisions be taken typically one-hour to one-day ahead of operations based on the uncertain data/information available at the time of decision.\par
The short-term operation of microgrids is fraught with a wide range of sources of uncertainty, including microgrid net load, i.e., microgrid load less \textit{VER} generation. In a scenario where a microgrid transacts energy on wholesale markets, the short-term planning may be exacerbated by the uncertainty associated with electricity market prices. As such, the judicious short-term planning of microgrids with integrated \textit{VER}s and potential exposure to the volatility in electricity market prices hinges on \textit{UC} approaches that undertake an explicit assessment of the uncertainty in net load and market prices.\par
To engage with uncertainty, most studies in the literature rely on stochastic optimization (\textit{SO}) or robust optimization (\textit{RO}) techniques. A major shortcoming of \textit{SO} is the assumption that the underlying probability distribution of uncertain parameters is known \textit{a priori}. The veracity of this assumption, however, is highly questionable, as system operators have access to collected data---not to their underlying probability distribution. Indeed, if the assumed probability distribution is incorrect, \textit{SO} may give rise to a markedly poor out-of-sample performance, which warrants and calls for optimization approaches that are not confined to a pre-specified probability distribution. In contrast to \textit{SO}, \textit{RO} techniques completely disregard the probabilistic nature of uncertain parameters and take decisions based solely on the worst-case scenario, which may yield overly conservative optimal solutions. \par
Distributionally robust optimization (\textit{DRO})---albeit being initially proposed long ago---has recently gained traction as a paradigm that addresses the drawbacks of both \textit{SO} and \textit{RO}. Under the \textit{DRO} paradigm, the probability distribution of uncertain parameters itself is considered to be uncertain and belong to an \textit{ambiguity set} of probability distributions that may be constructed based on various methods, including using moment information \cite{drogmb:45, drogmb:39}, the Kullback-Leibler (\textit{KL}) divergence \cite{drogmb:29}, and Wasserstein distance \cite{drogmb:46}. Central to \textit{DRO} is the formulation of an optimization problem that minimizes the expected cost brought about by the worst-case distribution in the ambiguity set. As such, \textit{DRO} obviates the need to commit to one pre-specified probability distribution and hedges the optimal decisions against adopting a misrepresenting probability distribution. The focus of this paper is the development of a \textit{DRO} approach for microgrid \textit{UC}.
\begin{table}%[!hbtp]
\begin{minipage}[c][7.8cm][t]{\linewidth}
\begin{fminipage}{\textwidth}[1ex]
\footnotesize
\renewcommand{\arraystretch}{1.2}
\begin{tabular}{l l  l l }
\multicolumn{2}{l}{\textbf{Nomenclature}}  & & \\
$\mathscr{H}$/$h$  & set/index of simulation time periods	&$\boldsymbol{x}$ & vector of first-stage variables \\
$\mathscr{G}$/$g$ & set/index of thermal generation 	&&  comprising $u_g[h]$ and $v_g[h]$	\\
& resources (\textit{TGR}s)	& $\boldsymbol{y}$ & vector of second-stage variables\\
$[p_{g} ]^{m}$/$[p_{g} ]^{M}$ & minimum/maximum power output	&& comprising $p_g[h]$, $p_b[h]$, and $p_s[h]$\\
&  of \textit{TGR} $g$	&$(\Omega, \mathcal{F}, \mathcal{P})$ & probability space\\
$T^{\uparrow}_{g} $/$T^{\downarrow}_{g} $ & minimum uptime/downtime	&$\mathcal{P}_{o}$ & nominal probability distribution\\
&   of \textit{TGR} $g$ 	&$\pi^{\omega}_{\mathcal{P}_{o}}$ & probability assigned to scenario $\omega$ \\
$c_{g}^{p}$/$c_{g}^{u}$ & linear/fixed fuel cost term for \textit{TGR} $g$	&&by the distribution $\mathcal{P}$\\
$c_{g}^{v}$ & start-up cost of \textit{TGR} $g$	&$\mathscr{P}$ & ambiguity set of \\
$p_{g}$ & power generation of \textit{TGR} $g$ in hour $h$	&&probability distributions\\
$u_{g}[h]$/$v_{g}[h]$ & binary comitment status/start-up	&$\boldsymbol{\tilde{\xi}}$ & random matrix associated with  \\
& variable of \textit{TGR} $g$ in hour $h$ 	& & net load and market prices\\
$p_{b}[h]$& power purchased from the electricity 	&$\boldsymbol{\xi^{\omega}}=[\boldsymbol{\eta^{\omega}}, \boldsymbol{\lambda^{\omega}}]$ & realization $\omega$ of $\boldsymbol{\tilde{\xi}}$\\
& market in hour $h$	&$\boldsymbol{\eta^{\omega}} \in \mathbb{R}^{24}$ & net load values in realization $\omega$\\
$p_{s}[h]$& spilled power in hour $h$	&$\boldsymbol{\lambda^{\omega}} \in \mathbb{R}^{24}$ & market prices in realization $\omega$\\
$\rho$  & divergence tolerance\\	
\end{tabular}
\end{fminipage}
\end{minipage}
\end{table}
\subsection{Related Work}
There is a growing body of literature on the application of \textit{DRO} approaches in the \textit{UC} problem. In \cite{drogmb:29}, the authors assess the uncertainty associated with wind generation and leverage the \textit{KL} divergence to propose a \textit{DRO} model for \textit{UC}; nevertheless, they do not evaluate the out-of-sample performance of their approach. The work conducted in \cite{drogmb:45} makes use of the first and second moment information to construct an ambiguity set and takes into account the uncertainty in \textit{VER} generation in the proposed \textit{UC} model. However, \cite{drogmb:45} does not consider the uncertainty associated with electricity prices. While \cite{drogmb:46} harnesses the Wasserstein distance, \cite{drogmb:39} capitalizes on moment information so as to construct ambiguity sets and develop \textit{DRO} approaches for \textit{UC}, yet neither approach jointly evaluates the uncertainty associated with \textit{VER} generation and electricity prices.\par
\subsection{Contributions and Structure of the Paper}
The general contributions and novel aspects of this paper are as follows:
\begin{enumerate}
\item We develop a new \textit{DRO} approach for microgrid \textit{UC} using the \textit{KL} divergence. To the best of our knowledge, this is the first study that jointly evaluates the uncertainty associated with microgrid net load and electricity market prices under a \textit{DRO} approach. We conduct representative studies and demonstrate the effectiveness of the proposed approach on real-world data.
\item Our studies provide valuable insights into the influence of divergence tolerance, and hence the degree of conservatism, on the out-of-sample performance.
\item We present a methodology that leverages the k-means clustering algorithm and soft dynamic time warping (\textit{SDTW}) score in constructing the nominal probability distribution and its support. The presented methodology lends itself to the joint study of the uncertainty in net load and electricity market prices through multidimensional clusters without unduly exacerbating the computational burden. 
\item We provide a tractable reformulation of the developed \textit{DRO} problem and present a two-level decomposition method in conjunction with an iterative algorithm that enables its solution by off-the-shelf solvers. The presented algorithm is amenable to parallelization on the basis of scenarios. 
\end{enumerate}
This paper contains four additional sections. In Section \ref{2}, we develop the mathematical formulation of the proposed \textit{UC} approach and spell out our ambiguity set construction methodology. We present an iterative decomposition method in Section \ref{3} for the solution of the proposed optimization problem. We illustrate the capabilities and effectiveness of the proposed \textit{DRO} framework in Section \ref{4} using representative studies and discuss the results. We present our concluding remarks in Section \ref{5}. 
\section{Mathematical Formulation}\label{2}
We devote this section to working out the analytical underpinning of the proposed \textit{DRO} approach. We discretize the time-axis and adopt 1 hour as the smallest indecomposable unit of time and 24 hours as the scheduling horizon. We define the study period by the set $\mathscr{H} \coloneqq \{h \colon h=1,...,24\}$.
%, which is designed to serve as a tool to conduct the short-term planning of a microgrid with an integrated \textit{TGR} and a \textit{VER}. We consider that the microgrid is interfaced with the main grid and has the capability to purchase power from the wholesale electricity market so as to meet its net load, i.e., microgrid load less \textit{VER} generation.
\subsection{Problem Formulation}\label{2a}
The proposed \textit{DRO} approach explicitly represents the uncertainty associated with net load and wholesale electricity market prices over the study period. We define by $\boldsymbol{\tilde{\xi}}$ the random matrix on the probability space $(\Omega, \mathcal{F}, \mathcal{P})$, where $\Omega$ is a sample space,
%is a finite sample space with cardinality $\vert\Omega\vert=\mathcal{S}$,%
$\mathcal{F}$ is a set of subsets of $\Omega$ that is a $\sigma-$algebra, and $\mathcal{P}$ is a probability distribution on $\mathcal{F}$. The random matrix $\boldsymbol{\tilde{\xi}} \in \Xi \subset \mathbb{R}^{24\times2}$ denotes the uncertain net load and electricity price values over the study period, where $\Xi$ denotes the support of the probability distribution $\mathcal{P}$. We assume that $\mathcal{P}$ has a finite support taking $\mathcal{S}$ realizations that we equivalently refer to as scenarios, i.e., $\vert\Xi\vert = \mathcal{S} < \infty$. The construction of $\Xi$ is detailed in Section \ref{2b}. For each realization $\omega \in \Omega$ of $\boldsymbol{\tilde{\xi}}$, we write the relation $\boldsymbol{\xi^{\omega}} = [\boldsymbol{{\eta}^{\omega}}, \boldsymbol{{\lambda}^{\omega}}]$, where $\boldsymbol{{\eta}^{\omega}} \in \mathbb{R}^{24}$ and $\boldsymbol{{\lambda}^{\omega}} \in \mathbb{R}^{24}$ represent the net load and electricity price values over the 24 hours of the study period $\mathscr{H}$, respectively. We denote by ${\xi^{\omega}}[h] = [{{\eta}^{\omega}}[h], {{\lambda}^{\omega}}[h]]$ the row $h$ of $\boldsymbol{\xi^{\omega}}$, which represents the net load (${{\eta}^{\omega}}[h]$) and electricity price (${{\lambda}^{\omega}}[h]$) in hour $h \in \mathscr{H}$ for scenario $\omega$. We denote by $\pi^{\omega}_{\mathcal{P}}$ the probability assigned to the scenario $\omega$ by the probability distribution $\mathcal{P}$. \par
The proposed formulation is based on a two-stage decision mechanism that mimics the order in which \textit{UC} and \textit{ED} decisions are taken.
\begin{IEEEeqnarray}{ll}
\underset{u_g[h], v_g[h]}{\text{minimize}} &  \bigg\{\sum_{h \in \mathscr{H}}  \sum_{g \in \mathscr{G}} c^{v}_{g} v_{g}[h] + c^{u}_{g} u_{g}[h]  \nonumber\\
&\hspace{30pt}+ \,\underset{\mathcal{P} \in \mathscr{P}}{\text{maximize}} \;\mathbb{E}_{\mathcal{P}}\big[\mathcal{Q}(\boldsymbol{x},\boldsymbol{\tilde{\xi}})\big] \bigg\}, \label{obj}\\
%&\hspace{30pt}+ \,\underset{\mathcal{P} \in \mathscr{P}}{\text{maximize}} \;\sum_{\omega\in \Omega}  \pi^{\omega}_{\mathcal{P}}\mathcal{Q}(\boldsymbol{x},\boldsymbol{\tilde{\xi}}) \bigg\}, \label{obj}\\
\text{subject to} & \nonumber\\
&v_{g}[h] \geq u_{g}[h]-u_{g}[h-1], \, \forall g \in \mathscr{G}, \forall h \in \mathscr{H},\label{st}\\
&u_{g}[h]  -  u_{g}[h-1]  \leq  u_{g}[\nu],\, \forall \nu \in \mathbb{N}\: \text{such that} \nonumber\\
&\quad \quad  h \leq \nu \leq \min\{h-1+T^{\uparrow}_{g},24\},\forall g \in \mathscr{G}, \label{upt}  \\
&u_{g}[h-1]  -  u_{g}[h]  \leq  1-u_{g}[\nu],\, \forall \nu \in \mathbb{N}\: \text{such that}  \nonumber\\
&\quad \quad  h \leq \nu \leq \min\{h-1+T^{\downarrow}_{g},24\},\forall g \in \mathscr{G}, \label{dwt}  \\
& u_{g}[h], v_{g}[h] \in \{0,1\}, \, \forall g \in \mathscr{G}, \forall h \in \mathscr{H}. \label{bin}
\end{IEEEeqnarray}
The first-stage problem \eqref{obj}-\eqref{bin} seeks to determine the binary commitment ($u_g[h]$) and start-up ($v_g[h]$) variables of the \textit{TGR}s over the study period, while taking into account the minimum uptime \eqref{upt} and downtime \eqref{dwt} constraints of the \textit{TGR}s. We represent all first-stage decision variables by the vector $\boldsymbol{x}$, which comprises $u_g[h]$ and $v_g[h]$.
%=[u_1[1], \ldots, u_1[24], \ldots, u_{G}[24], v_1[1], \ldots, v_1[24], \ldots, v_{G}[24]]^{\mathsf{T}}$. 
The first-stage decisions are taken before the realization of the uncertain net load values and electricity prices with the objective \eqref{obj} to minimize the fixed generation and start-up costs plus the worst-case expected power generation and purchase costs. \par
A salient feature of the proposed \textit{DRO} approach is to capitalize on an ambiguity set of probability distributions denoted by $\mathscr{P}$ to study the uncertainty associated with net load and electricity prices. Such an approach ensures that all probability distributions that belong to the set $\mathscr{P}$ be assessed and the optimal first-stage decisions be taken based on the expected cost brought about by the worst-case distribution. We elaborate on the construction of the set $\mathscr{P}$ in Section \ref{2b}.\par
The function $\mathcal{Q}(\boldsymbol{x},\boldsymbol{\tilde{\xi}})$ in \eqref{obj} denotes the uncertain power generation and purchase costs. For a specific vector of first-stage decision variables $\boldsymbol{x}^{\dagger}$ and a realization $\boldsymbol{\xi^{\omega}}$, $\mathcal{Q}(\boldsymbol{x^{\dagger}},\boldsymbol{{\xi}^{\omega}})$ is evaluated by solving the following second-stage problem: \\
\begin{IEEEeqnarray}{cl}
\hspace{-15pt}\mathcal{Q}(\boldsymbol{x^{\dagger}},\boldsymbol{{\xi}^{\omega}})\, \coloneqq   &\nonumber\\
 \hspace{-15pt}\underset{p_g[h], p_b[h], p_s[h]}{\text{minimize}} &  \bigg\{\sum_{h \in \mathscr{H}}  \sum_{g \in \mathscr{G}} {c}^{p}_{g}p_{g}[h] + \lambda^{\omega}[h]p_{b}[h] \bigg\}, \label{obj2}\\
 \hspace{-15pt}\text{subject to} & \nonumber\\
&\hspace{-5pt}u_{g}^{\dagger}[h]p_{g}^{m} \leq p_{g}[h] \leq u_{g}^{\dagger}[h]p_{g}^{M}, \nonumber\\
&\hspace{95pt}\forall g \in \mathscr{G}, \forall h \in \mathscr{H}\label{ol1}\\
%&\hspace{-5pt} p_{g}[h] \leq u_{g}^{\dagger}[h]p_{g}^{M}, \nonumber\\
%&\hspace{95pt}\forall g \in \mathscr{G}, \forall h \in \mathscr{H}\label{ol2}\\
%&\hspace{-5pt}p_{g}[h]-p_{g}[h-1] \leq \Delta^{\uparrow}_{g}u_{g}[h-1]  \nonumber\\
%&+ p_{g}^{M}(1-u_{g}[h-1]), \, \forall g \in \mathscr{G}, \forall h \in \mathscr{H},\label{ru}\\
%&\hspace{-5pt}p_{g}[h-1]-p_{g}[h] \leq \Delta^{\downarrow}_{g}u_{g}[h]  \nonumber\\
%&  \hspace{20pt}+ p_{g}^{M}(1-u_{g}[h]),\, \forall g \in \mathscr{G}, \forall h \in \mathscr{H}\label{rd}\\
&\hspace{-5pt}\sum_{g \in \mathscr{G}} p_{g}[h] + p_{b}[h] - p_{s}[h] = \eta^{\omega}[h], \nonumber\\
&\hspace{123pt}\forall h \in \mathscr{H},\label{pb}\\
&\hspace{-5pt}p_b[h],  p_s[h] \geq 0, \, \forall h \in \mathscr{H}. \label{nn}
\end{IEEEeqnarray}
The second-stage problem \eqref{obj2}-\eqref{nn} seeks to 
%determine the power dispatch of \textit{TGR}s ($p_g[h]$), purchased power from the wholesale electricity market ($p_b[h]$) and the spilled power ($p_s[h]$) over the study period with the objective \eqref{obj2} to 
minimize the power generation and purchase costs while taking into account the \textit{TGR} output limits \eqref{ol1} and power balance constraint \eqref{pb}.
\subsection{Ambiguity Set Construction Methodology}\label{2b}
We devote this subsection to the description of the methodology undertaken in constructing $\mathscr{P}$. %We start out by the development of the nominal probability distribution $\mathcal{P}_{o}$. 
We denote by $\mathcal{N}_{D}$ the number of days for which historical net load and electricity market price data are initially considered. We leverage the k-means clustering algorithm to partition the $\mathcal{N}_{D}$ number of multidimensional time-series data points to $\mathcal{S}$ multidimensional clusters so as to assign each data point to the cluster with the nearest cluster centroid. To this end, we make use of the \textit{SDTW} score to measure the similarity between time-series data points, which---when applied jointly with the k-means algorithm---was reported to deliver better results for time-series clustering tasks vis-à-vis the Euclidean distance \cite{drogmb:106}.\par
We note that the computational complexity of \eqref{obj}-\eqref{bin} gets aggravated with increasing number of uncertain parameters and scenarios. As such, we specifically aim at the joint representation of uncertain net load and market prices by multidimensional clusters, which affords the capability to simultaneously assess the uncertainty associated with net load and electricity prices without undue computational burden.\par
We use each of the constructed $\mathcal{S}$ clusters to form each of the $\mathcal{S}$ scenarios of the nominal probability distribution. For each cluster $\omega$, we utilize the cluster centroid to represent the realization $\boldsymbol{\xi^{\omega}}$ and construct the support $\Xi\coloneqq\{\boldsymbol{\xi^{\omega}}\colon\omega=1,\ldots,\mathcal{S}\}$. We denote by $\mathcal{N}_{\omega}$ the number of data points assigned to cluster $\omega$ and---for the nominal probability distribution $\mathcal{P}_{o}$---assign the probability for the scenario $\omega$ as $\pi^{\omega}_{\mathcal{P}_{o}}=\frac{\mathcal{N}_{\omega}}{\mathcal{N}_{D}},\,\omega=1,\ldots,\mathcal{S}$. \par
We next leverage the \textit{KL} divergence to construct an ambiguity set of probability distributions $\mathscr{P}$ around the nominal probability distribution $\mathcal{P}_{o}$. The ambiguity set formulation using the \textit{KL} divergence \cite{drogmb:108} is stated as:
\begin{IEEEeqnarray}{ll}
\mathscr{P}\coloneqq\{\mathcal{P} : & \sum_{\omega=1}^{S}\pi^{\omega}_{\mathcal{P}}\log\bigg(\frac{\pi^{\omega}_{\mathcal{P}}}{\pi^{\omega}_{\mathcal{P}_{o}}}\bigg)\leq\rho, \label{dc}\\
& \sum_{\omega=1}^{S}\pi^{\omega}_{\mathcal{P}}=1, \label{sc}\\
& \pi^{\omega}_{\mathcal{P}}\geq0\quad\forall \omega \in \Omega\}. \label{pc}
\end{IEEEeqnarray}
The divergence tolerance $\rho$ is a cornerstone of ambiguity set construction, in that it adjusts the size, and thus the degree of conservatism, of an ambiguity set. When $\rho=0$, $\mathscr{P}$ shrinks to a singleton that contains only the nominal distribution $\mathcal{P}_{o}$. On the flip side, as $\rho \rightarrow \infty$, $\mathscr{P}$ admits all probability distributions, which may result in overly conservative decisions.\par
For notational brevity, we present the \textit{KL} divergence-based microgrid \textit{UC} ($\mathrm{KL-MUC}$) formulation as:
\begin{IEEEeqnarray}{ll}
 \hspace{-15pt}\mathrm{KL-MUC:}   &\nonumber\\
 \hspace{-15pt}\underset{\boldsymbol{x}}{\text{minimize}} & \hspace{7pt} \boldsymbol{c}\cdot\boldsymbol{x} +\underset{\mathcal{P} \in \mathscr{P}}{\text{maximize}} \;\sum_{\omega \in \Omega}\pi^{\omega}_{\mathcal{P}}\mathcal{Q}(\boldsymbol{x},\boldsymbol{{\xi}^{\omega}}),  \label{obj3}\\
 \hspace{-15pt}\text{subject to} & \hspace{7pt}x \in \mathscr{X} , \label{xc}\\
& \hspace{7pt} \eqref{dc}-\eqref{pc}, \nonumber
\end{IEEEeqnarray}
where $\mathscr{X}$ represents the feasibility region of $\boldsymbol{x}$ defined by the constraints \eqref{st}-\eqref{bin}.
\section{Solution Method}\label{3}
In this section, we present a method based on Benders' decomposition to ensure the efficient solution of the $\mathrm{KL-MUC}$ problem. We take the dual of the inner maximization problem in $\mathrm{KL-MUC}$ and assign the dual variables $\zeta$ and $\mu$ to the constraints \eqref{dc} and \eqref{sc}, respectively, which, as per \cite{drogmb:108}, yields the following convex mixed-integer nonlinear reformulated $\mathrm{KL-MUC}$ ($\mathrm{RKL-MUC}$) problem:
\begin{IEEEeqnarray}{ll}
 \hspace{-15pt}\mathrm{RKL-MUC}:  &\hspace{-3pt}\nonumber\\
 \hspace{-15pt}\underset{\boldsymbol{x}, \mu, \zeta}{\text{minimize}} & \hspace{-3pt} \boldsymbol{c}\cdot\boldsymbol{x} + \mu + \rho \zeta + \zeta\sum_{\omega =1}^{\mathcal{S}} \pi^{\omega}_{\mathcal{P}_{o}}e^{\overline{\mathcal{K}}^{\omega}(\boldsymbol{x}, \mu, \zeta)-1},  \label{obj4}\\
 \hspace{-15pt}\text{subject to} & \hspace{-3pt}x \in \mathscr{X} , \label{xc}\\
 \hspace{0pt}& \hspace{-3pt}\zeta \geq 0, \label{lc}
\end{IEEEeqnarray}
where $\overline{\mathcal{K}}^{\omega}(\boldsymbol{x}, \mu, \zeta)=\frac{\mathcal{Q}(\boldsymbol{x},\boldsymbol{\xi^{\omega}})-\mu}{\zeta}$. For notational brevity, we define the following functions: 
\begin{IEEEeqnarray}{ll}
{\overline{\mathcal{R}}^{\omega}}(\boldsymbol{x},\zeta, \mu)  & \coloneqq  \zeta e^{\overline{\mathcal{K}}^{\omega}(\boldsymbol{x}, \mu, \zeta)-1},\label{rbarfnc}\\
{\mathcal{R}}(\boldsymbol{x},\zeta, \mu) & \coloneqq \sum_{\omega =1}^{\mathcal{S}} \pi^{\omega}_{\mathcal{P}_{o}}{\overline{\mathcal{R}}^{\omega}}(\boldsymbol{x},\zeta, \mu).\label{rfnc}
%{\mathcal{R}}(\boldsymbol{x},\zeta, \mu) &  \coloneqq  \zeta\sum_{\omega =1}^{\mathcal{S}} \pi^{\omega}_{\mathcal{P}_{o}}e^{\overline{\mathcal{K}}^{\omega}(\boldsymbol{x}, \mu, \zeta)-1},\nonumber\\
%& \hspace{3pt}=  \sum_{\omega =1}^{\mathcal{S}} \pi^{\omega}_{\mathcal{P}_{o}}{\overline{\mathcal{R}}^{\omega}}(\boldsymbol{x},\zeta, \mu).\label{rfnc}
\end{IEEEeqnarray}
We decompose $\mathrm{RKL-MUC}$ to a lower-bounding master problem ($\mathrm{MP}$) and an upper-bounding subproblem ($\mathrm{SP}$). 
\vspace{-0.2cm}
\begin{IEEEeqnarray}{ll}
 \hspace{0pt}  \mathrm{MP:}   &\hspace{6pt}\nonumber\\
 \hspace{0pt} \underset{\boldsymbol{x_{(\nu)}}, \mu_{(\nu)}, \zeta_{(\nu)}}{\text{minimize}} &\hspace{8pt}  \boldsymbol{c}\cdot\boldsymbol{x_{(\nu)}} + \mu_{(\nu)} + \rho \zeta_{(\nu)} + \theta_{(\nu)} \label{obj5}\\
\hspace{0pt}  \text{subject to} &\hspace{6pt} x_{(\nu)} \in \mathscr{X} , \label{xc}\\
 \hspace{0pt} &\hspace{6pt} \zeta_{(\nu)} \geq 0, \label{lc}\\
 \hspace{0pt}  & \hspace{6pt}\theta_{(\nu)} \geq \boldsymbol{{\alpha}_{(j)}}\cdot  (\boldsymbol{x_{(\nu)}}-\boldsymbol{x_{(j)}}) + {\beta}_{(j)} (\mu_{(\nu)}- \mu_{(j)})\nonumber\\
\hspace{35pt} & \hspace{35pt}+  {\gamma}_{(j)} (\zeta_{(\nu)}-\zeta_{(j)}) + {\mathcal{R}}(\boldsymbol{x_{(j)}}, \zeta_{(j)}, \mu_{(j)}), \nonumber\\
 \hspace{35pt}&\hspace{35pt}\,j =1,\ldots,\nu-1, \label{optc}
\end{IEEEeqnarray}
where \eqref{optc} represents the Benders’ optimality cuts that serve to approximate from below the function ${\mathcal{R}}(\boldsymbol{x} ,\zeta, \mu)$. At each iteration $\nu$, the candidate optimal variables $(\boldsymbol{x_{(\nu)}}, {\mu_{(\nu)}}, {\zeta_{(\nu)}})$ evaluated by the $\mathrm{MP}$ are fixed as $\boldsymbol{x_{f}}\gets\boldsymbol{x_{(\nu)}},$
${\mu_{f}}\gets {\mu_{(\nu)}}$, and ${\zeta_{f}}\gets{\zeta_{(\nu)}}$.
%\begin{IEEEeqnarray}{l}
%\boldsymbol{x_{f}}\gets\boldsymbol{x_{(\nu)}},\nonumber\\
%{\mu_{f}}\gets {\mu_{(\nu)}},\nonumber\\
%{\lambda_{f}}\gets{\lambda_{(\nu)}}.\nonumber
%\end{IEEEeqnarray}
Note that ${\mathcal{R}}(\boldsymbol{x}, \zeta, \mu)$ is a nonlinear convex function. To ease the computational burden, instead of minimizing ${\mathcal{R}}(\cdot)$ in the $\mathrm{SP}$, we adopt the linear program $\mathcal{Q}(\cdot)$ presented in \eqref{obj2}-\eqref{nn} as the $\mathrm{SP}$ \cite{drogmb:108} and leverage the chain rule along with the optimal $\mathrm{SP}$ solution to evaluate the optimality cuts. At each iteration $\nu$, the $\mathrm{SP}$ for each scenario $\omega$ is defined as: 
\begin{IEEEeqnarray}{lll}
 \hspace{-15pt}\mathrm{SP:}   &\hspace{7pt} &\nonumber\\
 \hspace{-15pt}\underset{\boldsymbol{y^{\omega}_{(\nu)}}, \boldsymbol{\hat{x}^{\omega}_{(\nu)}}}{\text{minimize}} &  \hspace{7pt} \sum_{h \in \mathscr{H}}  \sum_{g \in \mathscr{G}} {c}^{p}_{g}p_{g}[h] + \lambda^{\omega}[h]p_{b}[h], \label{obj22}\\
 \hspace{-15pt} \text{subject to} & \hspace{7pt} \eqref{ol1}-\eqref{nn},\nonumber\\
 & \hspace{7pt} \boldsymbol{\hat{x}^{\omega}_{(\nu)}}= \boldsymbol{x_{f}}& \hspace{-50pt}\xleftrightarrow{}\hspace{20pt} \boldsymbol{\overline{\varphi}^{\omega}_{(\nu)}}.\label{alc}
 \end{IEEEeqnarray}
The dual variable $ \boldsymbol{\overline{\varphi}^{\omega}_{(\nu)}}$ associated with the constraint \eqref{alc} represents the negative of the sensitivity of \eqref{obj22} to $\boldsymbol{x_{f}}$. We remark that the $\mathrm{SP}$ for each scenario $\omega$ is a continuous problem as the elements of $\boldsymbol{\hat{x}^{\omega}_{(\nu)}}$ are not constrained to be binary. For each scenario $\omega$, we evaluate the terms:
\begin{IEEEeqnarray}{ll}
\boldsymbol{\overline{\alpha}^{\omega}_{(\nu)}}&=\zeta_{f} e^{\overline{\mathcal{K}}^{\omega}(\boldsymbol{x_{f}}, \mu_{f}, \zeta_{f})-1}\, \boldsymbol{\overline{\varphi}^{\omega}_{(\nu)}}, \label{alphaeq}\\
\overline{\beta}^{\omega}_{(\nu)}&=\frac{\partial {\overline{\mathcal{R}}^{\omega}}(\boldsymbol{x_{f}},\zeta_f, \mu_f)}{\partial \zeta_f}\nonumber\\
&=(1-\overline{\mathcal{K}}^{\omega}(\boldsymbol{x_{f}}, \mu_{f}, \zeta_{f}))e^{\overline{\mathcal{K}}^{\omega}(\boldsymbol{x_{f}}, \mu_{f}, \zeta_{f})-1},\label{betaeq}	\\
\overline{\gamma}^{\omega}_{(\nu)}&=\frac{\partial {\overline{\mathcal{R}}^{\omega}}(\boldsymbol{x_{f}}, \zeta_f, \mu_f)}{\partial \mu_f}=-e^{\overline{\mathcal{K}}^{\omega}(\boldsymbol{x_{f}}, \mu_{f}, \zeta_{f})-1},\label{gammaeq}
\end{IEEEeqnarray}
and compute the terms $\boldsymbol{{\alpha}_{(\nu)}}$, ${{\beta}_{(\nu)}}$, ${{\gamma}_{(\nu)}}$ in the Benders' optimality cuts in \eqref{optc} as follows: $\boldsymbol{{\alpha}_{(\nu)}}=\sum_{\omega =1}^{\mathcal{S}} \pi_{\mathcal{P}_{o}}^{\omega} \boldsymbol{\overline{\alpha}^{\omega}_{(\nu)}}$, ${{\beta}_{(\nu)}}=\sum_{\omega =1}^{\mathcal{S}} \pi_{\mathcal{P}_{o}}^{\omega} {\overline{\beta}^{\omega}_{(\nu)}}$, ${{\gamma}_{(\nu)}}=\sum_{\omega=1}^{\mathcal{S}} \pi_{\mathcal{P}_{o}}^{\omega} {\overline{\gamma}^{\omega}_{(\nu)}}$. The feasibility cuts are not required in the $\mathrm{MP}$, as the $\mathrm{KL-MUC}$ problem has a relatively complete recourse and the \textit{KL} divergence does not necessitate a feasibility cut.\par
The presented method lends itself to parallelization on the basis of scenarios, since the terms $\boldsymbol{\overline{\alpha}^{\omega}_{(\nu)}}$, ${\overline{\beta}^{\omega}_{(\nu)}}$, and $ {\overline{\gamma}^{\omega}_{(\nu)}}$ can be computed independently for each scenario. We succinctly represent the decomposition algorithm in Algorithm \ref{algo}.\par
\begin{algorithm}
\caption{Decomposition algorithm for $\mathrm{RKL-MUC}$}\label{algo}
  \begin{algorithmic}[1]
    \STATE \text{Initialize} $\boldsymbol{x}\gets\boldsymbol{0}$. 
    \STATE \text{Solve} $\mathrm{SP}$. \text{Set} $\mathcal{Q}^{M}\gets \max(\{\mathcal{Q}(\boldsymbol{x}, \boldsymbol{{\xi}^{\omega}})\colon \omega \in \Omega\})$
    \STATE \text{Initialize} $\mathrm{UB}\gets\infty$, $\mathrm{LB}\gets-\infty$, $\nu\gets1$, $\theta_{(1)}\gets0$, $\zeta_{(1)}\gets0$, $\boldsymbol{x_{(1)}}\gets\boldsymbol{0}$.
    \WHILE {$\mathrm{UB}-\mathrm{LB}\geq\mathrm{TOL}$}
    	\STATE \text{Solve} $\mathrm{MP}$. \text{Determine} $\boldsymbol{x_{(\nu)}}, \zeta_{(\nu)}, \mu_{(\nu)}$, and $\theta_{(\nu)}$ so that  $\frac{\mathcal{Q}^{M}-\mu_{(\nu)}}{\zeta_{(\nu)}}\leq\mathcal{K}^{M}$. $\mathrm{LB}\gets\theta_{(\nu)}$.
	\STATE \text{Solve} $\mathrm{SP}$. \text{Determine} $\boldsymbol{{\alpha}_{(\nu)}}$, ${\beta}_{(\nu)}$, ${\gamma}_{(\nu)}$, and ${\mathcal{R}}(\boldsymbol{x_{(\nu)}}, \zeta_{(\nu)}, \mu_{(\nu)})$. $\mathrm{UB}\gets {\mathcal{R}}(\boldsymbol{x_{(\nu)}}, \zeta_{(\nu)}, \mu_{(\nu)})$.  $\nu\gets \nu+1$. 
    \ENDWHILE
  \end{algorithmic}
\end{algorithm}
We point out that $\boldsymbol{\overline{\alpha}^{\omega}_{(\nu)}}$, $\overline{\beta}^{\omega}_{(\nu)}$, and $\overline{\gamma}^{\omega}_{(\nu)}$ contain the term $\overline{\mathcal{K}}^{\omega}(\cdot)$ in the exponent, which renders the proposed method prone to overflowing errors during its execution. As such, we expressly stipulate a computational upper bound on $\overline{\mathcal{K}}^{\omega}(\cdot)$ denoted by $\mathcal{K}^{M}$. Nevertheless, in lieu of relying on $\overline{\mathcal{K}}^{\omega}(\cdot) \leq \mathcal{K}^{M}$ to bound $\overline{\mathcal{K}}^{\omega}(\cdot)$, we impose a more restrictive upper bound, \textit{viz}.: $\frac{\mathcal{Q}^{M}-\mu}{\zeta}\leq\mathcal{K}^{M}$. In contrast to \cite{drogmb:108} that evaluates $\overline{\mathcal{K}}^{\omega}(\cdot)$ in the $\mathrm{SP}$ and requires additional iterations to compute a new $\mu$ in the event that $(\boldsymbol{x_{(\nu)}}, {\mu_{(\nu)}}, {\zeta_{(\nu)}})$ evaluated by the $\mathrm{MP}$ prompts $\overline{\mathcal{K}}^{\omega}(\cdot)$ to be greater than $\mathcal{K}^{M}$, our proposed upper bound ensures that $(\boldsymbol{x_{(\nu)}}, {\mu_{(\nu)}}, {\zeta_{(\nu)}})$ determined by the $\mathrm{MP}$ satisfy $\overline{\mathcal{K}}^{\omega}(\cdot)\leq \mathcal{K}^{M}$ $\forall \omega \in \Omega$, thereby precluding the need for additional iterations.
\section{Case Study and Results}\label{4}
In this section, we carry out representative studies to illustrate the application and effectiveness of the proposed \textit{DRO} approach. We consider a microgrid with an integrated \textit{TGR} and a \textit{PV} panel. The source code and simulation scripts for the case study are provided in \cite{drogmb:111}. The load and \textit{PV} generation dataset \cite{drogmb:80} contains measurements collected from June 1, 2019 to August 31, 2019 in an anonymous house in New York. To ensure consistency, we consider the locational marginal prices at the N.Y.C. bus in the New York Independent System Operator network cleared in the day-ahead market for the said time period and add a surcharge to the prices so as to reflect the rates available to residential customers \cite{drogmb:110}. \par
%We split the dataset so as to utilize the measurements in June and July to solve the proposed \textit{DRO} formulation with different values of divergence parameter $\rho$ as well as the corresponding \textit{SO} formulation---which serves as a benchmark for our experiments. Next, we fix the values of the first-stage decision variables for each experiment and utilize the data for August 2019 to evaluate the out-of-sample performance of each approach.\par
We start out by the construction of the scenarios. We utilize the data collected from June 1, 2019 to July 31, 2019 and deploy the methodology described in Section \ref{2b} to assign each data point to $\mathcal{S}$ clusters. To determine $\mathcal{S}$, we examine the
% influence of $\mathcal{S}$ by harnessing the Silhouette and Davies Bouldin metrics---all the while evaluating the 
percentage of variance captured for different values of $\mathcal{S}$ and pick $\mathcal{S}=8$, at which the point of diminishing returns (i.e., the so-called \textit{elbow}) is reached and 88.06\% of the total variance is captured, where capturing an additional 10\% of the variance requires 36 more clusters.\par
%We employ the formed clusters to construct the nominal probability distribution and the scenarios for our experiments.\par
We draw on the solution method described in Section \ref{3} to solve the $\mathrm{RKL-MUC}$ problem. We perform our implementations in \textit{Pyomo} using Gurobi 9.0.2 as the solver with the optimality tolerance gap $\mathrm{TOL}=10^{-5}$ on a 2.6 \textit{GHz} Intel Core i7 \textit{CPU} with 16 \textit{GB} of \textit{RAM}. The discussion in Section \ref{2b} revealed the divergence tolerance $\rho$ as a key determinant of the degree of conservatism of $\mathscr{P}$. As such, we probe the influence of $\rho$ by solving the $\mathrm{RKL-MUC}$ problem with each of the following values: $\rho=0, 0.2, 0.4, 0.6, 0.8,1.0$. To carry out comparative assessments, we develop the following equivalent stochastic formulation of the $\mathrm{KL-MUC}$ problem, which serves as a benchmark for our experiments: 
\begin{IEEEeqnarray}{ll}
\hspace{-15pt}\mathrm{SUC:}  &\nonumber\\
\hspace{-15pt}\underset{\boldsymbol{x}}{\text{minimize}} & \hspace{7pt} \boldsymbol{c}\cdot\boldsymbol{x} +\sum_{\omega \in \Omega}\pi^{\omega}_{\mathcal{P}_{o}}\mathcal{Q}(\boldsymbol{x},\boldsymbol{{\xi}^{\omega}}),  \label{obj7}\\
\hspace{-15pt}\text{subject to} & \hspace{7pt}x \in \mathscr{X} , \label{xc2}
\end{IEEEeqnarray}
and employ the L-shaped algorithm for its solution.\par
While the $\mathrm{RKL-MUC}$ and $\mathrm{SUC}$ problems are solved using the constructed ambiguity sets, their feasibility must be assessed on real-life data that were not harnessed in constructing the ambiguity sets. To this end, %we empirically investigate the out-of-sample performance of $\mathrm{RKL-MUC}$ and $\mathrm{SUC}$. Out-of-sample performance, in the strictest sense of the term, entails to an analytical evaluation using the true underlying data-generating distribution. Since the data-generating distribution is not known, we empirically investigate the out-of-sample performance using real-life data. 
we capitalize on the data collected from August 1, 2019 to August 31, 2019 to form the out-of-sample dataset and empirically investigate the out-of-sample performance of the $\mathrm{RKL-MUC}$ formulation for each of the considered six values of $\rho$, as well as that of the $\mathrm{SUC}$ problem. To do so, for each of the seven setups, we fix the optimal first-stage decisions obtained using the constructed ambiguity sets and scenarios, and we subsequently compute the total cost by providing each setup with the data points of the out-of-sample dataset. \par
\begin{figure}[!h]
\vspace{0cm}
\centering
\begin{tikzpicture}
\pgfplotsset{}
\begin{axis}[
xlabel={\textit{divergence tolerance }$\rho$},
ylabel={\textit{total cost (\textit{\$})}},
width=0.64\textwidth,
height=0.40\textwidth,
xlabel style={at={(axis description cs:0.5, 0.03)},anchor=north},
ylabel style={at={(axis description cs:0.03,.5)},rotate=0, anchor=south},
xmin=0,
xmax=1.0,
xtick = {0,0.2,0.4,0.6, 0.8, 1.0},
ymin=1320,
ymax=1400,
ytick = {1320,1340, 1360, 1380, 1400},
legend style={font=\scriptsize, at={(0.22,0.35)}, legend columns=1, anchor=north},
label style={font=\small},
tick label style={font=\small}  ]
	\addplot+[sharp plot, white!10!olive, solid, semithick, mark=o] plot coordinates
		{ (0, 1380.7691121120004) (0.2, 1380.7691121120004) (0.4, 1380.7691121120004) (0.6, 1380.7691121120004) (0.8, 1380.7691121120004) (1.0, 1380.7691121120004)  };
	\addlegendentry{$\mathrm{SUC}$ \hspace{0.12cm}}
	\addplot+[sharp plot, black!30!red, solid, semithick, mark=triangle] plot coordinates
		{ (0, 1380.7691121120004) (0.2, 1380.7691) (0.4, 1361.5880070790004) (0.6, 1345.1975793000004) (0.8, 1345.6935524510004) (1.0, 1345.6935524510004)  };
	\addlegendentry{$\mathrm{RKL-MUC}$\hspace{0.12cm}}
	\end{axis}
\end{tikzpicture}
\vspace{0cm}
\caption{Out-of-sample performances under $\mathrm{RKL-MUC}$ and $\mathrm{SUC}$}
\vspace{0cm}
\label{res}
\end{figure}
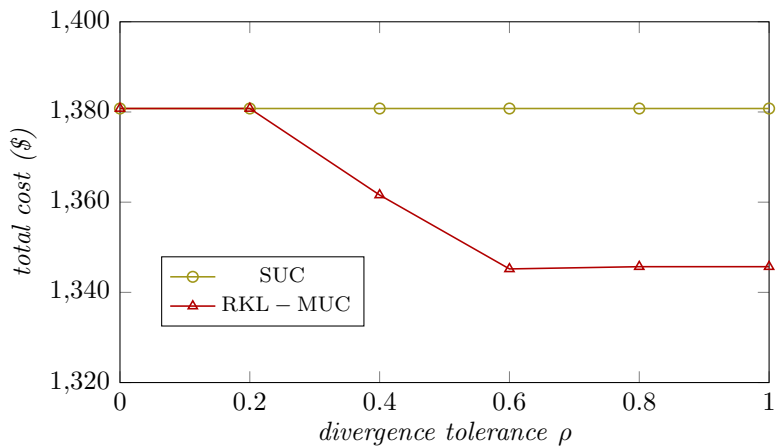

We present in Fig. \ref{res} the total cost under the $\mathrm{RKL-MUC}$ and $\mathrm{SUC}$ formulations. 
%While the solid lines in Fig. \ref{res} indicate the mean value of the daily total costs of each day in the out-of-sample dataset, the shaded regions around the solid lines demarcate one standard deviation around the mean value. \par
At the outset, we note that the $\mathrm{SUC}$ solution tallies with the $\mathrm{RKL-MUC}$ solution for $\rho=0$, which validates our computations, as when $\rho=0$, the ambiguity set contains solely the nominal probability distribution and so the $\mathrm{KL-MUC}$ formulation reduces to the $\mathrm{SUC}$ formulation. We remark upon the fact that, for all considered $\rho$ values, the total cost under the $\mathrm{RKL-MUC}$ formulation is less than or equal to that under the $\mathrm{SUC}$ formulation. We further observe that the total cost decreases as $\rho$ increases from $0$ to $0.6$. 
%In fact, as $\rho$ increases from $0$ to $0.6$, the standard deviation of the daily total costs analogously decreases, which indicates that the daily total costs become more concentrated around the mean value and less spread out. 
These observations bring out the benefit of taking into account additional probability distributions other than the nominal probability distribution and make clear that the nominal probability distribution need not be taken at face value. This notwithstanding, the total cost slightly picks up as $\rho$ increases above $0.6$, which may be accounted for by the fact that the assignment of increasingly large values to $\rho$ permits the incorporation of probability distributions that assign markedly high probabilities to adverse scenarios into the ambiguity set, which are evidently not reflected in the out-of-sample dataset.

%\begin{figure}[!h]
%\vspace{0cm}
%\centering
%\begin{tikzpicture}
%\pgfplotsset{}
%\begin{axis}[
%xlabel={\textit{divergence tolerance }$\rho$},
%ylabel={\textit{total costs (\textit{\$})}},
%width=0.46\textwidth,
%height=0.26\textwidth,
%xlabel style={at={(axis description cs:0.5, 0.03)},anchor=north},
%ylabel style={at={(axis description cs:0.03,.5)},rotate=0, anchor=south},
%xmin=0,
%xmax=1.0,
%xtick = {0,0.2,0.4,0.6, 0.8, 1.0},
%ymin=1342,
%ymax=1382,
%ytick = {1342,1352, 1362, 1372, 1382},
%legend style={font=\scriptsize, at={(0.22,0.35)}, legend columns=1, anchor=north},
%label style={font=\small},
%tick label style={font=\small}  ]
%	\addplot+[white!10!olive, smooth, solid, semithick, mark=o] plot coordinates
%		{ (0, 1380.7691121120004) (0.2, 1380.7691121120004) (0.4, 1380.7691121120004) (0.6, 1380.7691121120004) (0.8, 1380.7691121120004) (1.0, 1380.7691121120004)  };
%	\addlegendentry{$\mathrm{SUC}$ \hspace{0.12cm}}
%	\addplot+[black!30!red, smooth, solid, semithick, mark=triangle] plot coordinates
%		{ (0, 1380.7691121120004) (0.2, 1380.7691) (0.4, 1361.5880070790004) (0.6, 1345.1975793000004) (0.8, 1345.6935524510004) (1.0, 1345.6935524510004)  };
%	\addlegendentry{$\mathrm{RKL-MUC}$\hspace{0.12cm}}
%	\end{axis}
%\end{tikzpicture}
%\vspace{-0.3cm}
%\caption{Out-of-sample performances under $\mathrm{RKL-MUC}$ and $\mathrm{SUC}$}
%\vspace{-0.2cm}
%\label{res}
%\end{figure}
\section{Conclusion}\label{5}
In this paper, we propose a \textit{DRO} approach for microgrid unit commitment under net load and electricity price uncertainty. Our approach takes full advantage of the copious amounts of data imparted by the deployment of information and communication technologies as per the smart grid paradigm. The methodology leveraged in constructing the scenarios affords the capability to conjointly study the uncertainty associated with net load and electricity prices without aggravating the computational burden. The hallmark of our approach is to minimize the worst-case expected cost over an ambiguity set of probability distributions constructed using the \textit{KL}-divergence, which enables us to hedge the optimal decisions against adopting a misrepresenting probability distribution. The case studies conducted on real-world data demonstrate the effectiveness of the proposed approach.

\end{document}